\documentclass{article}
\usepackage{geometry}[top=15truemm, bottom=15truemm, right=10truemm, left=10truemm]
\usepackage{amsmath,amssymb,amsthm}
\usepackage{amsmath,amssymb,amsthm}
\usepackage[all]{xy}
\usepackage{graphicx}
\usepackage{mathrsfs}
\usepackage{lineno}

\newcommand{\pp}{\mathcal{P}}
\newcommand{\ee}{\mathcal{E}}

\newcommand{\idp}{\mathfrak{p}}
\newcommand{\idq}{\mathfrak{q}}

\newcommand{\Ad}{\mathbb{A}}

\newcommand{\nat}{\mathbb{N}}

\theoremstyle{plain}
\newtheorem{thm}{Theorem}[section]
\newtheorem{lem}[thm]{Lemma}
\newtheorem{prop}[thm]{Proposition}
\newtheorem{cor}[thm]{Corollary}

\theoremstyle{definition}
\newtheorem{defn}[thm]{Definition}

\newtheorem{rem}[thm]{Remark}

\numberwithin{equation}{section}
\title{Absolute reconstruction of number fields\\ from the Deligne-Ribet monoids}
\author{Takeo Uramoto\\ Graduate School of Science and Engineering, Kagoshima University}
\date{}
\begin{document}
\maketitle
\begin{abstract}
\noindent
Following Cornelissen, Li, Marcolli, and Smit, this short paper proves that the field structure of a number field $K$ can be reconstructed from the pair $(DR_K, I_K)$ of the Deligne-Ribet monoid $DR_K$ and the submonoid $I_K \subseteq DR_K$ when $K$ is the rational number field or an imaginary quadratic field. The general-case reconstruction is also discussed, which is more abstract than the case of rational and imaginary quadratic fields. 
\end{abstract}
\section{Introduction}
\label{s1}
Cornelissen, Li, Marcolli, and Smit \cite{Cornelissen_Li_Marcolli_Smit} proved that two number fields $K, L$ are isomorphic as fields if and only if there exists a topological-monoid isomorphism $\phi: DR_K \rightarrow DR_L$ of their Deligne-Ribet monoids in such a way that $\phi$ restricts to the monoid isomorphism $\phi: I_K \rightarrow I_L$ of the submonoids $I_K \subseteq DR_K$, $I_L \subseteq DR_L$ of the non-zero integral ideals of $K, L$. The current paper then describes how to reconstruct the field structure of $K$ from the semigroup structure of $DR_K$ and $I_K \subseteq DR_K$, where $I_K$ plays an auxiliary role to relate the first and second factors of $DR_K = \hat{O}_K \times_{\hat{O}_K^\times} G_K^{ab}$; this reconstruction also implies yet another proof of the main result of \cite{Cornelissen_Li_Marcolli_Smit} at least for the rational and imaginary quadratic fields; the general-case reconstruction is also developed, which however heavily relies on the results of Cornelissen et.\ al.\ \cite{Cornelissen_Li_Marcolli_Smit} and Hoshi \cite{Hoshi} and is far more abstract than the case of rational and imaginary quadratic fields (cf.\ \S \ref{s4}). 

Our method to reconstruct the field $K$ from $DR_K$ is based on some elementary ideas from the theory of profinite semigroups, which will highlight how the semigroup structure of $DR_K$ is useful to reconstruct the arithmetic structure of the number field $K$. A basis for our method is to classify the idempotents of $DR_K$ in terms of sets of the maximal ideals $\idp \in P_K$; we will see that certain idempotents $e_S$ of $DR_K$, which are explicitly constructed from $S \subseteq P_K$ and played a key role in \cite{Cornelissen_Li_Marcolli_Smit} too, actually \emph{exhaust} all idempotents of $DR_K$. That is:

\begin{prop}
Let $S \subseteq P_K$ be any subset of the set $P_K$ of maximal ideals of the integer ring $O_K$ of $K$; and define an idempotent $e_S \in DR_K$ as follows (cf.\ \S 7 \cite{Cornelissen_Li_Marcolli_Smit}, or \S \ref{s2} below): 
\begin{eqnarray}
 e_S &:=& [1_S, 1],
\end{eqnarray}
Then every idempotent $e \in DR_K$ is equal to $e_{S_e}$ for some unique subset $S_e \subseteq P_K$, which is explicitly given as follows: writing $e = [\rho_e, s_e] \in \hat{O}_K \times_{\hat{O}_K^\times} G_K^{ab} \simeq DR_K$, 
\begin{eqnarray}
 S_e &:=& \bigl \{ \idp \in P_K \mid v_\idp(\rho_{e, \idp}) = 0 \bigr\}.
\end{eqnarray}
\end{prop}

In other words, we have a bijective correspondence between the set $\ee_K$ of idempotents of $DR_K$ and the set $\pp(P_K)$ of the subsets of $P_K$. We also note that a natural poset-structure on idempotents of $DR_K$ can be defined by $e \leq e'$ if and only if $e \cdot e' = e$, which is in fact equivalent to the inclusion $S_e \subseteq S_{e'}$. Therefore, the above correspondence $\ee_K \rightarrow \pp(P_K)$ is actually an isomorphism of posets; in particular, by the poset structure, it makes sense to define an idempotent $e \in \ee_K$ to be \emph{maximal} by saying that (i) $e\neq 1$ and (ii) $e \leq e'$ implies that $e'=e$ or $e'=1$. 

With this concept, we can classify the maximal ideals $\idp \in P_K$ of $O_K$ in the following form:

\begin{cor}
There exists a bijective correspondence between the maximal ideals $\idp \in P_K$ of $O_K$ and the maximal idempotents of $DR_K$. 
\end{cor}

Starting from this classification, we discuss how to reconstruct the field structure of the number field $K$ purely from the semigroup structure of $DR_K$. Our next key step is to observe that the multiplicative monoid $O_\idp$ of the integers in the local field $K_\idp$ is naturally embedded as a submonoid into the first factor of $DR_K = \hat{O}_K\times_{\hat{O}_K^\times} G_K^{ab}$, whose image can be characterized semigroup-theoretically as follows:

\begin{prop}
An element $x \in DR_K$ belongs in the image of the natural embedding $O_\idp \hookrightarrow DR_K$ if and only if the following identities hold: 
\begin{eqnarray}
 x \cdot e_\emptyset &=& e_\emptyset\\
 x \cdot e_{\{ \idq \}}&=& e_{\{ \idq\}}
\end{eqnarray}
for all $\idq \neq \idp$. 
\end{prop}
\noindent
We can also give similar characterizations for the multiplicative monoids $O_\idp^\times$ and $O_\idp^* := O_\idp \setminus \{0\}$, and thus, the multiplicative group $K_\idp^\times$ of the local field $K_\idp$ as well, since $K_\idp^\times$ is the groupification of $O_\idp^*$; notice further that, for each $\rho_\idp \in O_\idp^*$, there exists a \emph{unique} element $\sigma_\idp \in DR_K^\times$ such that $\rho_\idp \cdot \sigma_\idp \in I_K$, with which we can define a homomorphism $O_\idp^* \ni \rho_\idp \mapsto \sigma_\idp \in DR_K^\times$. Then, combining these observations, we obtain the following diagram:
\begin{equation}
\xymatrix{
O_\idp^* \ar[r] \ar[d]_{\textrm{grp.}} & (\Ad_K^\times/ \overline{K^\times \cdot K^{\times, \circ}_\infty}) \\
K^\times_\idp \ar@{-->}[ur] & 
}
\end{equation}
where the dotted map $K_\idp^\times \rightarrow (\Ad_K^\times/ \overline{K^\times \cdot K^{\times, \circ}_\infty})$ coincides with the usual reciprocity map. 
Therefore, by reconstructing the group structure of the finite idele group $\Ad_{K, f}^\times$ from these $K_\idp^\times$, we obtain the following description of the group $K^\times$ (recall that $K$ is supposed to be an imaginary quadratic or the rational number field, whence $\overline{K^\times \cdot K^{\times, \circ}_\infty} = K^\times \cdot K^{\times, \circ}_\infty$): For an imaginary quadratic field $K$, 
\begin{eqnarray}
\label{why imaginary quadratic}
 K^\times &=& \ker (\Ad_{K, f}^\times \rightarrow (\Ad_K^\times/ \overline{K^\times \cdot K^{\times, \circ}_\infty}));
\end{eqnarray}
(The case of the rational number field is similar, except for some care about the archimedian part; we shall only deal with imaginary quadratic fields in this paper.) From this, the local ring $O_{K, \idp} \subseteq K$ too can be reconstructed by the fact that $O_{K, \idp}^* = O_{K, \idp}\setminus \{0\}$ is equal to the intersection of $K^\times$ and $O_\idp^*$ within $K_\idp^\times$. 

In this way, we can reconstruct the field $K = K^\times \cup \{0\}$ from $DR_K$ at least as a monoid, together with additional information on the monoid structure of the local rings $O_{K, \idp}$; we then apply the result of Hoshi \cite{Hoshi} to reconstruct the field structure of $K$ from these data. As a consequence of this, we can also deduce that a topological-monoid isomorphism $\phi: DR_K \rightarrow DR_L$ gives a field isomorphism $\phi: K \rightarrow L$ in a canonical way. 

The general-case reconstruction is also discussed in the last section (\S \ref{s4}). This method is, however, more abstract than the case of the rational and imaginary quadratic fields. Some more explicit method that works for general case should be pursued further; we shall left this problem unsolved in this paper.

\paragraph{Acknowledgements}
We are grateful to Masanori Morishita for his continuous encouragements. This work was originally prepared for some part of our talk given at Kagoshima University, 28, July 2025; we are grateful to the organizers and colleagues for providing us this valuable opportunity; the result in \S \ref{s4} is our answer to a question from Kunio Obitsu there, for which we are also grateful to him. This work is supported by JSPS KAKENHI No.22K03248. 
\section{The local structure of $DR_K$}
\label{s2}
This section proves some preparatory results on the semigroup structure of $DR_K$. To this end, let us start with recalling some basic constructions of idempotents in profinite monoids. Throughout this paper, for a number field $K$, we shall denote by $P_K$ the set of maximal ideals of $K$ and by $I_K$ the monoid of non-zero integral ideals of $K$. For the definition of the Deligne-Ribet monoids, the reader is referred to \cite{Cornelissen_Li_Marcolli_Smit} or \cite{Deligne_Ribet}; also for some concepts and facts in semigroup theory, we refer the reader to \cite{Rhodes_Steinberg}. 

\begin{defn}[cf.\ \cite{Rhodes_Steinberg}]
For a profinite monoid $M$ and an arbitrary element $a \in M$, the sequence $(a^{n!})_{n \in \nat}$ in $M$ converges to an idempotent in $M$ which is denoted by $a^\omega$: That is,
\begin{eqnarray}
a^\omega &:=& \lim_{n \to \infty} a^{n!}.
\end{eqnarray}
\end{defn}

\begin{rem}[cf.\ \cite{Rhodes_Steinberg}]
It is well-known that $a \in M$ is invertible if and only if $a^\omega = 1$, which we shall use below. In particular, for any $s \in G_K^{ab}$, we have $[1, s]^\omega = 1$ in $DR_K$. 
\end{rem}

\begin{defn}[normal form of idempotents]
Let $S \subseteq P_K$ be any subset of $P_K$. Then define the idempotent $e_S \in \ee_K$ as follows:
\begin{eqnarray}
 e_S &:=& [1_S, 1],
\end{eqnarray}
where $1_S \in \hat{O}_K$ is such that $(1_S)_\idp = 1$ if $\idp \in S$ and $(1_S)_\idp = 0$ otherwise. 
\end{defn}

\begin{prop}[classification of idempotents]
\label{classification of idempotents}
For an idempotent $e \in \ee_K$, let $S_e \in \pp(P_K)$ be defined as follows: denoting $e = [\rho_e, s_e] \in \hat{O}_K \times_{\hat{O}_K^\times} G_K^{ab}$, 
\begin{eqnarray}
 S_e &:=& \bigl\{ \idp \in P_K \mid v_\idp (\rho_{e, \idp}) = 0\}. 
\end{eqnarray}
Then the correspondences $e \mapsto S_e$ and $S \mapsto e_S$ are mutually inverse. 
\end{prop}
\begin{proof}
The equality $S=S_{e_S}$ holds since $v_\idp (0) = \infty > 0$ and $v_\idp (1) = 0$. Conversely, we shall prove the equality $e = e_{S_e}$ for any idempotent $e \in DR_K$. Denote $e = [\rho_e, s_e] \in \hat{O}_K \times_{\hat{O}_K^\times} G_K^{ab} \simeq DR_K$ for some $\rho_e \in \hat{O}_K$ and $s_e \in G_K^{ab}$. Since $e = e^\omega$, we have:
\begin{eqnarray}
 e &=& [\rho_e, s_e]^\omega \\
 &=& [\rho_e, 1]^\omega \cdot [1, s_e]^\omega \\
&=& [\rho_e, 1]^\omega.
\end{eqnarray}
Notice here that the $\idp$-component $[\rho_{e, \idp}, 1]^{n!}$ converges to $0$ if and only if $v_\idp (\rho_{e, \idp}) > 0$, that is, if and only if $\idp \not \in S$; otherwise, $[\rho_{e, \idp}, 1]^{n!}$ converges to $1$ because $\rho_{e, \idp} \in O_\idp^\times$ then. This proves that $e = e_{S_e}$ as requested. 
\end{proof}

\begin{lem}
For any $S, S' \subseteq P_K$, we have:
\begin{eqnarray}
 e_S \cdot e_{S'}  &=& e_{S \cap S'}.
\end{eqnarray}
In particular, $e_S \leq e_{S'}$ if and only if $S \subseteq S'$. 
\end{lem}
\begin{proof}
The first claim is clear from the definition of $e_S$, from which the second claim also readily follows. 
\end{proof}

\begin{cor}[maximal idempotents and maximal ideals]
\label{maximal idempotents and maximal ideals}
An idempotent $e \in DR_K$ is maximal if and only if the complement $P_K \setminus S_e$ of the corresponding $S_e$ within $P_K$ is a singleton. Therefore the following map gives a bijection $e \mapsto \idp_e$ from the maximal idempotents of $DR_K$ to the maximal ideals of $O_K$:
\begin{eqnarray}
 e &\longmapsto& P_K \setminus S_e = \{\idp_e\}. 
\end{eqnarray}
\end{cor}

\section{Reconstruction of number field}
\label{s3}
\noindent
Based on the above results on the semigroup structure of $DR_K$, we develop a semigroup-theoretic method to reconstruct the number field $K$ from $DR_K$. For this purpose, we start with proving a few lemmas:

\begin{lem}
An element $x \in DR_K$ is representable as $x = [\rho, 1]$ with some $\rho \in \hat{O}_K$ if and only if we have the following identity:
\begin{eqnarray}
 x \cdot e_{\emptyset} &=& e_\emptyset. 
\end{eqnarray}
\end{lem}
\begin{proof}
The only-if part is trivial; to prove the if part, suppose that we have $x = [\rho, s]$ with $\rho \in \hat{O}_K$, $s \in G_K^{ab}$, and the identity $x \cdot e_{\emptyset}=e_\emptyset$. Then, since $e_\emptyset = [0, 1]$, we have $[0, s] = [0, 1]$, which means that there exists some $u \in \hat{O}_K^\times$ such that $s = u$; thus, we deduce $x = [\rho, s] = [\rho u, 1]$. This completes the proof. 
\end{proof}

\begin{lem}
An element $x \in DR_K$ belongs in the image of the embedding $O_\idp \hookrightarrow DR_K$ if and only if we have the following identities:
\begin{eqnarray}
 x \cdot e_{\emptyset} &=& e_\emptyset \\
 x \cdot e_{\{\idq\}} &=&  e_{\{\idq\}},
\end{eqnarray}
for any $\idq \neq \idp$. 
\end{lem}
\begin{proof}
This follows from the above lemma and the definition of the idempotent $e_{\{\idq\}}$. 
\end{proof}

\begin{cor}
An element $x \in DR_K$ belongs in the image of $O_\idp^\times$ under the embedding $O_\idp \hookrightarrow DR_K$ if and only if $x \in O_\idp$ and $x^\omega=1$. 
\end{cor}

\begin{cor}
An element $x \in DR_K$ belongs in the image of $O_\idp^*$ under the embedding $O_\idp \hookrightarrow DR_K$ if and only if $x \in O_\idp$ and $x \neq 0_\idp$, where $0_\idp \in O_\idp$ is the unique element of $O_\idp$ such that $y \cdot 0_\idp = 0_\idp$ for all $y \in O_\idp$. 
\end{cor}

Now we construct a canonical monoid homomorphism $O^*_\idp \rightarrow DR_K^\times$, where the submonoid $I_K \subseteq DR_K$ plays a key role: 

\begin{lem}
For each $\rho_\idp \in O^*_\idp$, there exists a unique $\sigma_\idp \in DR_K^\times$ such that $\rho_\idp \cdot \sigma_\idp \in I_K$. 
\end{lem}
\begin{proof}
For $\rho_\idp \in O^*_\idp$, define $\sigma_\idp = [1, \rho_\idp^{-1}] \in DR_K^\times$, whence we have $\rho_\idp \cdot \sigma_\idp = [\rho_\idp, \rho_\idp^{-1}] \in I_K$. To see the uniqueness, suppose that $\sigma = [1, s] \in DR_K^\times$ is such that $[\rho_\idp, s] \in I_K$. Then, there exists some $\rho \in \hat{O}_K^* \cap \Ad_{K, f}^\times$ such that $[\rho_\idp, s] = [\rho, \rho^{-1}]$, that is, there exists a unit $u \in \hat{O}_K^\times$ such that $\rho_\idp = \rho u$ and $[s] = [(\rho u)^{-1}] \in G_K^{ab}$, thus, $[s] = [\rho_\idp^{-1}] = \sigma_\idp$ in particular. This completes the proof. 
\end{proof}

Therefore, we can define a monoid homomorphism $O^*_\idp \rightarrow DR_K^\times$ by $\rho_\idp \mapsto \sigma_\idp$, which then induces a group homomorphism $K_\idp^\times \rightarrow DR_K^\times$ since $K_\idp^\times$ is now characterized as the groupification of $O_\idp^*$. Using this, we can reconstruct the multiplicative group $K^\times$ as follows:

\begin{prop}
Define $\Ad_{K, f}^\times$  from the groupifications $K_\idp^\times$ of $O_\idp^*$ in the usual way. Then the unit group $K^\times$ of the number field $K$ is isomorphic to the kernel of the following homomorphism defined naturally from the above maps $K_\idp^\times \rightarrow DR_K^\times =  (\Ad_K^\times/ \overline{K^\times \cdot K^{\times, \circ}_\infty})$: 
\begin{equation}
\label{finite reciprocity map}
\xymatrix{
\Ad_{K, f}^\times \ar[r] &  (\Ad_K^\times/ \overline{K^\times \cdot K^{\times, \circ}_\infty})
}\end{equation}
\end{prop}
\begin{proof}
In our restricted case, we have $\overline{K^\times \cdot K^{\times, \circ}_\infty} = K^\times \cdot K^{\times, \circ}_\infty$. Thus the kernel of (\ref{finite reciprocity map}) is isomorphic to the multiplicative group $K^\times$. 
\end{proof}

Now our main result is the following: 

\begin{thm}[reconstruction of field structure]
The field structure of $K$ can be reconstructed from the semigroup structure of $DR_K$. 
\end{thm}
\begin{proof}
We need to reconstruct the additive structure of the multiplicative monoid $K = K^\times \cup \{0\}$; to this end, notice first that we have not only the multiplicative monoid $K$ but also its submonoids $O_{K, \idp} \subseteq K$ for each $\idp \in P_K$ by noting the following characterization of $O_{K, \idp}^* = O_{K, \idp}\setminus \{0\}$:
\begin{equation}
\xymatrix{
O_{K, \idp}^* \ar[r] \ar[d] & O_\idp^* \ar[r] \ar[d] & (\Ad_K^\times/ \overline{K^\times \cdot K^{\times, \circ}_\infty}) \\
K^\times \ar[rd] \ar[r] & \Ad_{K, f}^\times \ar[ur] \ar@{->>}[d] & \\
& K_\idp^\times
}
\end{equation}
where $O_{K, \idp}^*$ is the intersection of $K^\times$ and $O_\idp$ within $K_\idp^\times$; from these, we also have $O_K = \bigcap_\idp O_{K, \idp}$. Second, recall that we can characterize the group $U_\idp^{(1)} = 1 + \idp O_\idp \leq O_\idp^\times$ as $U_\idp^{(1)} = (O_\idp^\times)^{N\idp -1}$, cf.\ \cite{Neukirch}, and thus, the group $O_{K, \idp}^{\prec} = 1 + \idp O_{K, \idp} \subseteq O_{K, \idp}^\times$ too. Hence, with these constructions, we eventually reconstruct the data $(K, O_K, P_K, \{ O_{K, \idp}^{\prec} \})_\idp$ from the semigroup structure of $DR_K$; using it we then reconstruct the additive (and thus field) structure of $K$ by the result of Hoshi \cite{Hoshi}. This completes the proof. 
\end{proof}

It follows from our constructions that a topological-monoid isomorphism $\phi: DR_K \rightarrow DR_L$ that restricts to an isomorphism $\phi: I_K \rightarrow I_L$ induces a field isomorphism $K \rightarrow L$ in a canonical way. 

\section{Concluding remarks}
\label{s4}
\noindent
After we introduced the above results in our talk at Kagoshima University, Kunio Obitsu asked us whether the absolute Galois group $G_K$ can be reconstructed from the Deligne-Ribet monoid $DR_K$. Our answer was as follows: Since we could reconstruct $K$ from $DR_K$ at least when $K$ is the rational or imaginary quadratic field, we can of course reconstruct from $DR_K$ the absolute Galois group $G_K$ too at least for such fields. But we noticed later that we might answer to his question affirmatively (yet somewhat trivially) for arbitrary number fields in such a way that heavily relies on the results of Cornelissen et.\ al.\ \cite{Cornelissen_Li_Marcolli_Smit} and Hoshi \cite{Hoshi}. 

Let $K$ be an arbitrary number field and $DR_K$ its Deligne-Ribet monoid. To reconstruct the field structure of $K$ from the semigroup structure of $DR_K$, it suffices to reconstruct a profinite group $G$ which is isomorphic to the absolute Galois group $G_K$ thanks to the work of Hoshi \cite{Hoshi} (probably as Obitsu intended by his question). In fact, we may say that this is possible by combining the results of Hoshi \cite{Hoshi} \emph{and} Cornelissen et.\ al.\ \cite{Cornelissen_Li_Marcolli_Smit}: Consider the set of profinite groups $G$ of \emph{AGSC-type} (cf.\ \S 3 \cite{Hoshi}) such that the profinite monoid $D(G) := \hat{O}(G) \times_{\hat{O}(G)^\times} G^{ab}$ together with the submonoid $I(G)$ is isomorphic to $DR_K$ as topological monoids so that $I(G) \simeq I_K$, where note that all the necessary data (e.g.\ suitable monoids $\hat{O}(G)$ and $I(G)$) can be constructed from $G$ \emph{purely group-theoretically}: Indeed, for a group-theoretic construction of the monoid $\hat{O}(G)$, see Theorem 1.4 (6) \cite{Hoshi}; also, for the submonoid $I(G)$, use Frobenius elements (cf.\ Theorem 1.4 (4) \cite{Hoshi}) for maximal closed subgroups of $G$ of MLF-type (cf.\ Proposition 3.5 \cite{Hoshi}), recalling e.g.\ $\idp = [\pi_\idp, \pi_\idp^{-1}]$ in $DR_K$ for $\idp \in P_K$. The point is that such a profinite group $G$ is \emph{uniquely determined} up to isomorphism by the result of \cite{Cornelissen_Li_Marcolli_Smit}. It is in this sense that we say that a profinite group $G$ isomorphic to $G_K$, thus the field structure of the number field $K$ too, can be reconstructed from the semigroup structure of $DR_K$ for arbitrary number field $K$. 

This procedure is more abstract than the explicit method in \S \ref{s2}--\S \ref{s3} in that, while the method in \S \ref{s2}--\S \ref{s3} is intrinsic to $DR_K$, the above one refers to the outer set of profinite groups $G$ of AGSC-type with $(D(G), I(G)) \simeq (DR_K, I_K)$; but one should note carefully that it concerns only abstract profinite groups $G$ and topological monoid isomorphisms $D(G) \simeq DR_K$. Thus, if this transcendental method (together with \cite{Hoshi}) is counted as a reconstruction of the field structure of a number field $K$ from the semigroup structure of $DR_K$, we are done; the use of the result of Cornelissen et.\ al.\ \cite{Cornelissen_Li_Marcolli_Smit} for this transcendental procedure would be allowed in view of the fact that Hoshi's method \cite{Hoshi} too relies on Neukirch-Uchida theorem. That said, it would be satisfactory if we could find some more explicit method in the general case as in the case of the rational and imaginary quadratic fields; this problem is not solved in this paper. 

\bibliographystyle{abbrv}
\bibliography{completeness_of_DR}
\end{document}